\input amstex
\documentstyle{amsppt}
\magnification=1200
\NoBlackBoxes
\hsize=5.8in
\vsize=7.5in
\topmatter
\title Conditional Association and Spin Systems\endtitle
\rightheadtext {Positive Correlations}
\date July 4, 2005\enddate
\author Thomas M. Liggett\endauthor
\thanks Research supported in part by
NSF Grant DMS-03-01795.\endthanks
\affil University of
California, Los Angeles\endaffil
\abstract A 1977 theorem of T. Harris states that an attractive spin system preserves
the class of associated probability measures. We study
analogues of this result for measures that satisfy various conditional positive
correlations properties. In particular, we show that a spin system
preserves measures satisfying the FKG lattice condition (essentially) if and 
only if distinct spins flip independently. The downward FKG property,
which has been useful recently in the study of the contact process,
lies between the properties of lattice FKG and association.
We prove that this property is preserved by a spin system if the death
rates are constant and the birth rates are additive (e.g., the contact
process), and prove a partial converse to this statement. Finally, we introduce
a new property, which we call downward conditional 
association, which lies between the FKG lattice condition and downward FKG, and
find essentially necessary and sufficient conditions for this property to
be preserved by a spin system. This suggests that the latter property may be
more natural than the downward FKG property.
\endabstract 
\endtopmatter

\heading 1. Introduction\endheading
Correlation inequalities have been used frequently in probability theory
and statistical physics. In this paper, we will consider correlation inequalities for
probability measures on $\{0,1\}^S$, where $S$ is a finite set, and
their connection with certain continuous time Markov chains
on $\{0,1\}^S$ that are known as spin systems. The key
definition is the following: the probability measure $\mu$ has
positive correlations, or is associated, if 
$$\int fgd\mu\geq\int fd\mu \int gd\mu\tag 1.1$$
for all increasing functions $f$ and $g$ on
$\{0,1\}^S$.
Two results proved in the 1970's provide convenient ways to show that a
probability measure on $\{0,1\}^S$ is associated:

\proclaim {Theorem 1.2} (FKG Theorem) Suppose $\mu$ is a probability measure
on $\{0,1\}^S$ that assigns strictly positive mass to every point in
$\{0,1\}^S$, and satisfies
$$\mu(\eta\wedge\zeta)\mu(\eta\vee\zeta)\geq\mu(\eta)\mu(\zeta)\tag 1.3$$
for all $\eta,\zeta\in\{0,1\}^S$. Then $\mu$ is associated.
\endproclaim
Assumption (1.3) is called the FKG lattice condition, or the strong FKG condition.

For the statement of the second of these results, recall that a spin system is a continuous time
Markov chain $\eta_t$ on
$\{0,1\}^S$ in which transitions can occur at only one site at a time. Let $\beta(x,\eta)$
and
$\delta(x,\eta)$ be the rates at which the transitions $0\rightarrow 1$ (births) and $1\rightarrow 0$
(deaths) occur at site $x$ if the configuration is $\eta$. (The
functions $\beta(x,\eta)$ and $\delta(x,\eta)$ do
not depend on $\eta(x)$.) The spin system is said to be attractive if
$\beta(x,\eta)$ is an increasing function of $\eta$ and $\delta(x,\eta)$ is a decreasing
function of $\eta$ for each $x$. Let $S(t)$ be the semigroup for the spin system.
It acts on functions and measures in the following way:
$$S(t)f(\eta)=E^{\eta}f(\eta_t);\qquad \int fd[\mu S(t)]=\int S(t)fd\mu.$$

\proclaim {Theorem 1.4} (Harris) If the spin system is attractive, then $\mu S(t)$ is
associated whenever $\mu$ is.
\endproclaim

See pages 78 and 80 of Liggett (1985) for proofs of Theorems 1.2 and 1.4. The second
of these implies that the stationary measure of an irreducible
attractive spin system is associated. Such a measure may or
may not satisfy (1.3). (It does satisfy (1.3) if the spin system is reversible
and attractive, but this is a very special situation.) A major advantage of
Theorem 1.4 over Theorem 1.2 is that one does not need an explicit expression
for the stationary distribution in order to check that it is associated.

More recently, some conditional forms of positive correlation inequalities have
been proved (Belitsky, Ferrari, Konno and Liggett (1997) and
van den Berg, H\"aggstr\"om and Kahn (2005a)) and applied (Liggett and
Steif (2005)). We are primarily concerned here with conditional versions of
Theorem 1.4. First, we state a converse to Theorem 1.4. It will be proved
in Section 2.

\proclaim{Theorem 1.5} Suppose a spin system has the property
that $\mu S(t)$ is associated whenever $\mu$ is. Then
the spin system is attractive.\endproclaim

An easy consequence of Theorem 1.2 is that (1.3) is equivalent to 
the property that not only $\mu$, but also all measures obtained
from $\mu$ by conditioning on the values of $\eta$ at any set of
sites in $S$, are associated. (We assume here and in the
sequel that all measures considered assign strictly positive probability
to the events on which one is conditioning.) Thus a first natural question is, for which spin systems is it
the case that $\mu S(t)$ satisfies (1.3) whenever $\mu$ does? The next
result says that this
occurs essentially only when the coordinate processes $\{\eta_t(x), x\in S\}$
are independent. Therefore, we see that
one should not condition on too much information if one hopes to have preservation
under the semigroup of some conditional positive correlations property for interesting
spin systems.
Say that the spin system has independent flips if for each $x$,
 $\beta(x,\eta)$ and
$\delta(x,\eta)$ do not depend on $\eta$.

\proclaim {Theorem 1.6} (a) Suppose that the spin system has independent flips.
Then $\mu S(t)$ satisfies (1.3) whenever $\mu$ does.

(b) Suppose that $\mu S(t)$ satisfies (1.3) whenever $\mu$ does, and that
$S$ has at least four points. Then the spin system has independent flips.
\endproclaim
Theorem 1.6 will be proved in Section 3. The statement in part (b) is trivially
true if $S$ is a singleton, but is false if $S$ has either two or three points.
If $S$ has two points, (1.3) is equivalent to association, so
that by Theorems 1.4 and 1.5, property (1.3) is preserved by the semigroup
if and only if the spin system is attractive -- it need not have
independent flips. If $S$ has three points, the following is an example
of a spin system that preserves (1.3) without having independent flips:
$\beta(x,\eta)=\delta(x,\eta)=0$ for all $\eta$, except that
$$\beta(x,\eta)=1\text{ if }\eta\equiv 1\quad\text{and}\quad \delta(x,\eta)
=1\text{ if }\eta\equiv 0.$$
To check this, note that $\mu S(t)(\eta)=e^{-t}\mu(\eta)$, except when $\eta\equiv 0$ or
$\eta\equiv 1$, and that $\mu S(t)(\eta)$ is increasing in $t$ if $\eta\equiv 0$ or $\eta\equiv 1$.
Of course, this example is reducible. To construct an irreducible example, simply
add a constant to all the rates. Since the process with constant rates preserves (1.3)
by Theorem 1.6(a), the modified process preserves (1.3) as well, by Proposition 1.14 below.

Following Liggett and Steif (2005), we will say that $\mu$ is downward FKG if
for any $A\subset S$, the conditional measure $\mu\{\cdot\mid\eta\equiv 0\text{ on }A\}$
is associated; i.e., conditioning is allowed only on
0's, not on 1's. This property lies between (1.3), where
one is allowed to condition on any configuration on $A$, and 
association, where one is not allowed to condition at all. Van den Berg,
H\"aggstr\"om and Kahn (2005a) proved that the distribution of the contact process at time
$t$ is downward FKG provided that the initial distribution is deterministic. 
(When we refer to such properties of measures on $\{0,1\}^S$ for
infinite $S$, we mean that all projections on finite subsets of $S$
have the property.) This property
was used by Liggett and Steif (2005) to show that the stationary distribution of the contact
process on $Z$ (and as a consequence, on many other graphs) dominates a nontrivial product measure. (The contact
process stationary distribution does not satisfy the FKG lattice condition -- see Liggett (1994)
and van de Berg, H\"aggstr\"om and Kahn (2005b), where it is shown that even conditioning
on $\eta(x)=1$ at a single site $x$ destroys the property of association.)

Another example of a measure that is downward FKG but does not satisfy the FKG lattice
condition is the following. Let $\pi$ be a permutation of $\{1,...,n\}$ that is chosen uniformly
at random. Define $\eta(i)$ to be the indicator of the event that $\pi(i)\neq i$, and
let $\mu$ be the distribution of $\eta$. Fishburn, Doyle and Shepp (1988) showed
that $\mu$ is associated, even though it does not satisfy the FKG lattice condition.
Since $\mu\{\cdot\mid \eta\equiv 0\text{ on }A\}$ when $A$ is of size $k$ corresponds to the
measure $\mu$ for random permutations of $n-k$ points, it follows that $\mu$
is downward FKG.

Our next result is an analogue
of Theorems 1.4 and 1.5 on the one hand, and Theorem 1.6 on the other, for preservation of
the downward FKG property. It will be proved in Section 4. 
For its statement, we need the following
definition. The birth rates $\beta(x,\eta)$ are said to
be additive if they can be written in the form
$$\beta(x,\eta)=\sum_{A\subset S}c(x,A)1_{\{\eta\not\equiv 0\text{ on } A\}}\tag 1.7$$
with $c(x,A)\geq 0$ for all $x,A$. The contact process of course has
additive birth rates and constant death rates.

\proclaim {Theorem 1.8} (a) Suppose that the spin system satisfies
$$\delta(x,\eta)\text{ does not depend on $\eta$ for each $x$}\tag 1.9$$
and the birth rates $\beta(x,\eta)$ are additive.
Then $\mu S(t)$ is downward FKG whenever $\mu$ is.

(b) Suppose that $\mu S(t)$ is downward FKG whenever $\mu$ is. Then 
$\delta(x,\eta)$ is constant on the set $\{\eta:\eta\not\equiv 0\}$ for each $x$.

(c) Suppose that all transition rates are zero, except
$\beta(u,\eta)$ for a particular $u\in S$ and arbitrary $\eta$. If $\mu S(t)$ is downward FKG whenever $\mu$
is, then $\beta(u,\eta)$ is increasing in $\eta$, and satisfies
$$\beta(u,\eta\vee\zeta)+\beta(u,\eta\wedge\zeta)\leq\beta(u,\eta)+\beta(u,\zeta)\tag 1.10$$
for all $\eta$ and $\zeta$.
\endproclaim

Following completion of the present paper, we learned that in the revision of
van de Berg, H\"aggstr\"om and Kahn (2005a), the authors added a result (Theorem 3.5)
that is part (a) of Theorem 1.8, restricted to the contact process. Our result
applies to a broader class of spin systems, and we believe that our
proof is significantly simpler. 

Note that (1.10) is an additive form of (1.3), and is satisfied whenever
$\beta(u,\cdot)$ is additive. However, there is a large gap between (1.10) and
additivity -- we do not have a necessary and sufficient condition for preservation
of the downward FKG property. Part (c) of the theorem does say that
monotonicity of the birth rates alone is not sufficient.

The fact that there is significant difference between our necessary conditions and sufficient conditions
for preservation of the downward FKG property suggests that this property may not be the
most natural one to consider.  We therefore introduce a new concept. The probability measure
$\mu$ will be called downward conditionally associated (DCA) if for every strictly positive
{\it decreasing} function $h$ on $\{0,1\}^S$ that satisfies the lattice condition
$$h(\eta\vee\zeta)h(\eta\wedge\zeta)\geq h(\eta)h(\zeta),\tag 1.11$$
the measure
$$\mu_h(d\eta)=h(\eta)\mu(d\eta)\bigg/\int hd\mu$$
is associated.  The implications among these properties are now:

$$\text{ FKG lattice}\quad\Rightarrow\quad \text{DCA}
\quad\Rightarrow\quad \text{downward FKG}\quad\Rightarrow\quad\text{association.}\tag 1.12$$
To check the first, note that as a result of (1.11), $\mu_h$ satisfies the
FKG lattice condition whenever $\mu$ does. For the second, given $A\subset S$,
apply the association of $\mu_h$ for 
$$h(\eta)=\prod_{x\in A}\big[1+\epsilon-\eta(x)\big],$$
and pass to the limit as $\epsilon\downarrow 0$. 

One should of course ask whether
the implications in (1.12) are strict. If $S$ has two points, it is easy to check that
they are all equivalences, and are equivalent to $\mu(11)\mu(00)\geq\mu(10)\mu(01)$.
In Section 6, we will give necessary and sufficient conditions for each of the four
properties appearing in (1.12) when $S$ has three points. A consequence of this
is that the first and third implications are strict in this case, while the
middle one is an equivalence. We do not know whether the middle implication is
strict for larger $S$. While this is an interesting question, it is in a sense not
too important from the point of view of this paper. Suppose, for example, that one
wanted to show that the contact process invariant measure is DCA. One approach is
to use Theorem 1.13 below. Another would be to use Theorem 1.8, and then show that
downward FKG implies DCA, if this is the case. Given the difficulty of showing this
implication when $S$ has three points, it seems clear that one should use Theorem 1.13 instead.
\medskip
\noindent {\bf Remark.} The monotonicity assumption on $h$ in the definition of DCA
is very important. To see this, consider the fact (which was pointed out
to the author by L. Chayes) that if $\mu_h$ is associated for
every positive $h$ that satisfies (1.11), then $\mu$ satisfies the FKG lattice condition.
To see this, suppose it does not. Then one can condition $\mu$ on the event $\{\eta:\eta\equiv \zeta\text{ on }
A\}$ for some $A\subset S$ and some $\zeta\in\{0,1\}^A$ and obtain a measure that
is not associated. It follows that $\mu_h$ is not associated if $h$ is some positive
perturbation of $1_{\{\eta:\eta\equiv\zeta\text{ on }A\}}$ satisfying (1.11).
\medskip
We are now in a position to state the following analogue of
Theorem 1.8 for the DCA property, but with conditions that are essentially necessary and sufficient.
Its proof will be given in Section 5.

\proclaim {Theorem 1.13} (a) Suppose that the death rates of the spin system satisfy
(1.9) and the birth rates are increasing and satisfy (1.10).
Then $\mu S(t)$ is DCA whenever $\mu$ is.

(b) Suppose that $\mu S(t)$ is DCA whenever $\mu$ is. Then 
$\delta(x,\eta)$ is constant on the set $\{\eta:\eta\not\equiv 0\}$ for each $x$.

(c) Suppose that all transition rates are zero, except
$\beta(u,\eta)$ for a particular $u\in S$. If $\mu S(t)$ is DCA whenever $\mu$
is, then $\beta(u,\eta)$ is increasing in $\eta$, and satisfies
(1.10) for all $\eta$ and $\zeta$.
\endproclaim

One consequence of Theorem 1.13 is that the upper invariant measure of the contact 
process has the presumably stronger DCA property. Another advantage over Theorem 1.8 is that
it can be used to prove this property for the stationary distributions of
 a larger class of spin systems.

We conclude the introduction by stating a result that will be used to simplify the proofs of Theorems
1.8 and 1.13.

\proclaim{Proposition 1.14} Suppose $\Cal L_1, \Cal L_2, S_1(t)$ and $S_2(t)$ are the generators
and semigroups for two Markov chains on $\{0,1\}^S$, and that $S(t)$ is the semigroup
with generator $\Cal L_1+\Cal L_2$. If $\Cal C$ is a closed set of functions or measures on
$\{0,1\}^S$ and $S_i(t)$ maps $\Cal C$ into itself for $i=1,2$, then $S(t)$ also
maps $\Cal C$ into itself.
\endproclaim

\demo{Proof} This is an immediate consequence of the Trotter product formula (Ethier and Kurtz (1986), page 33):
$$S(t)=\lim_{n\rightarrow\infty}\big[S_1(t/n)S_2(t/n)\big]^n.$$
\enddemo

\heading 2. A converse to Harris' Theorem\endheading
In this section we prove a stronger form of Theorem 1.5 which will also
be useful in later sections. It implies Theorem 1.5 because every product measure
is associated. This statement is an immediate consequence of either
Theorem 1.2 or Theorem 1.4. We will use the following notation:
If $\eta\in\{0,1\}^S$ and $x_1,x_2,...$ are distinct elements of $S$, then $\eta_{x_1,x_2,...}$
is the configuration with
$$\eta_{x_1,x_2,...}(x)=\cases1-\eta(x)&\text{ if }x=x_i\text{ for some }i\\
\eta(x)&\text{ otherwise.}\endcases$$

\proclaim{Proposition 2.1} Suppose the semigroup for a spin system has the property
that $\mu S(t)$ is associated for every product
measure $\mu$. Then the spin system is attractive.
\endproclaim
\demo{Proof}
Fix distinct $x,y\in S$, and let 
$\gamma$ be any configuration with $\gamma(x)=\gamma(y)=0$.
Let $\mu_{\epsilon}$ be the product measure with marginals
$$\mu_{\epsilon}\{\eta:\eta(z)=1\}=\cases\rho\quad&\text{if }z=x,\\
\lambda\quad&\text{if }z=y,\\\epsilon\quad&\text{if }z\neq
x,y,\gamma(z)=0\\1-\epsilon\quad&\text{if }z\neq x,y,\gamma(z)=1,
\endcases$$
where $0<\rho,\lambda,\epsilon<1$. 
By assumption, $\mu_{\epsilon}S(t)$ is associated for all $t$. Applying the
definition of association to the increasing functions $f(\eta)=\eta(x)$
and $g(\eta)=\eta(y)$ gives
$$\gathered\mu_{\epsilon}S(t)\{\eta(x)=1, \eta(y)=1\}\mu_{\epsilon}S(t)\{\eta(x)=0,
\eta(y)=0\}\\-\mu_{\epsilon}S(t)\{\eta(x)=1, \eta(y)=0\}\mu_{\epsilon}S(t)\{\eta(x)=0,
\eta(y)=1\}\geq 0
\endgathered\tag 2.2$$
Since the left side of (2.2) is zero at $t=0$, its derivative is nonnegative at
$t=0$. Differentiating (2.2) with respect to $t$, setting $t=0$, and then letting
$\epsilon\downarrow 0$ leads to the following inequality:
$$\gather\rho\lambda\bigg[\rho(1-\lambda)\delta(x,\gamma)+(1-\rho)\lambda\delta(y,\gamma)
-(1-\rho)(1-\lambda)\beta(x,\gamma)-(1-\rho)(1-\gamma)\beta(y,\gamma)\bigg]\\
+(1-\rho)(1-\lambda)\bigg[(1-\rho)\lambda\beta(x,\gamma_y)+\rho(1-\lambda)
\beta(y,\gamma_x)-\rho\lambda\delta(x,\gamma_y)-\rho\lambda\delta(y,\gamma_x)\bigg]\\
\geq\rho(1-\lambda)\bigg[(1-\rho)(1-\lambda)\beta(y,\gamma)+\rho\lambda\delta(x,\gamma_y)
-(1-\rho)\lambda\beta(x,\gamma_y)-(1-\rho)\lambda\delta(y,\gamma)\bigg]\\
+(1-\rho)\lambda\bigg[(1-\rho)(1-\lambda)\beta(x,\gamma)+\rho\lambda\delta(y,\gamma_x)
-\rho(1-\lambda)\delta(x,\gamma)-\rho(1-\lambda)\beta(y,\gamma_x)\bigg].
\endgather$$
Dividing by $\rho(1-\rho)\lambda(1-\lambda)$ and collecting like terms gives:
$$\frac{\beta(x,\gamma_y)-\beta(x,\gamma)}{\rho}+\frac{\beta(y,\gamma_x)-
\beta(y,\gamma)}{\lambda}+\frac{\delta(x,\gamma)-\delta(x,\gamma_y)}{1-\rho}
+\frac{\delta(y,\gamma)-\delta(y,\gamma_x)}{1-\lambda}\geq 0.$$
Now let respectively $\rho\rightarrow 0,\lambda\rightarrow 0,\rho\rightarrow 1,
\lambda\rightarrow 1$ to conclude that
$$\beta(x,\gamma)\leq\beta(x,\gamma_y),\quad\beta(y,\gamma)\leq\beta(y,\gamma_x),\quad
\delta(x,\gamma_y)\leq\delta(x,\gamma),\quad\delta(y,\gamma_x)\leq\delta(y,\gamma).$$
Since this is true for all $x,y,\gamma$ satisfying $\gamma(x)=\gamma(y)=0,$
it follows that the spin system is attractive.
\enddemo

\heading 3. Preservation of the FKG lattice condition\endheading
In this section, we prove Theorem 1.6. It will be convenient to recall that
(1.3) holds for all $\eta,\zeta$ if and only if it holds whenever
$\eta$ and $\zeta$ differ at exactly two sites.
\demo{Proof of Theorem 1.6(a)} Let $P_t(\eta,\gamma)$ be the transition
probabilities for the spin system. Then
$$\mu S(t)(\gamma)=\sum_{\eta}\mu(\eta)P_t(\eta,\gamma).\tag 3.1$$
For each $z\in S$, let $p_t(z,0,0),p_t(z,0,1),p_t(z,1,0),p_t(z,1,1)$ be the transition
probabilities for the two state Markov chain that has transitions
$0\rightarrow 1$ and $1\rightarrow 0$ at rates $\beta(z,\eta)$ and
$\delta(z,\eta)$ respectively. (Recall that we are assuming that these rates do not depend on
$\eta$.) Then
$$P_t(\eta,\gamma)=\prod_{z\in S}p_t(z,\eta(z),\gamma(z)).\tag 3.2$$
Fix two distinct sites $x,y\in S$ and let $\gamma$ be a configuration that
satisfies $\gamma(x)=\gamma(y)=0$. We must show that if
$\mu$ satisfies (1.3), then
$$\mu S(t)(\gamma_{x,y})\mu S(t)(\gamma)- \mu S(t)(\gamma_x)\gamma S(t)(\gamma_y)\geq 0.
\tag 3.3$$
Using (3.1), the left side of (3.3) can be written as
$$\aligned \frac 12\sum_{\eta,\zeta}\mu(\eta)\mu(\zeta)\bigg[&P_t(\eta,\gamma_{x,y})
P_t(\zeta,\gamma)+P_t(\zeta,\gamma_{x,y}) P_t(\eta,\gamma)\\&-P_t(\eta,\gamma_{x})
P_t(\zeta,\gamma_y)-P_t(\zeta,\gamma_{x}) P_t(\eta,\gamma_y)\bigg].
\endaligned\tag 3.4$$
Using (3.2), the expression in brackets in (3.4) becomes (after some cancellation)
$$\prod_{z\neq x,y}p_t(z,\eta(z),\gamma(z))p_t(z,\zeta(z),\gamma(z))\big[f(\eta)-f(\zeta)
\big]\big[g(\eta)-g(\zeta)\big],$$
where 
$$f(\eta)=p_t(x,\eta(x),1)\quad\text{and}\quad g(\eta)=p_t(y,\eta(y),1).$$
Note that $f$ and $g$ are increasing functions, so that by Theorem 1.2 applied
to the measure $\nu$, where
$$\nu(\eta)=c\mu(\eta)\prod_{z\neq x,y}p_t(z,\eta(z),\gamma(z)),$$
(3.4) is nonnegative. (Here $c$ is a normalizing constant.) The measure $\nu$ satisfies
(1.3) since $\mu$ does, and
$$\frac{\nu(\eta\wedge\zeta)\nu(\eta\vee\zeta)}{\nu(\eta)\nu(\zeta)}=
\frac{\mu(\eta\wedge\zeta)\mu(\eta\vee\zeta)}{\mu(\eta)\mu(\zeta)}.$$
\enddemo

We will isolate the main part of the proof of Theorem 1.6(b) in the following
proposition, since it will also be useful in Sections 4 and 5.
\proclaim{Proposition 3.5} Suppose that $\mu S(t)$ is downward FKG whenever
$\mu$ satisfies (1.3). Then $\delta(x,\eta)$ is constant on the set
$\{\eta:\eta\not\equiv 0\}$ for each $x\in S$.
\endproclaim
\demo{Proof} Any product measure satisfies (1.3), and downward FKG implies association,
so we can apply Proposition 2.1 to conclude that $\delta(x,\eta)$
is decreasing in $\eta$. 

Now take three distinct sites $x,y,z$ and let $\mu_{\epsilon}$ be the probability
measure on $\{0,1\}^S$ with respect to which $\{\eta(w),w\in S\backslash\{y,z\}\}$ are
independent with 
$\mu_{\epsilon}\{\eta:\eta(w)=1\}=\frac 12$
and independently of these, $(\eta(y),\eta(z))$ takes the following values:
$$\align(1,1)&\quad\text{with probability } 1-3\epsilon\\
(0,1)&\quad\text{with probability } \epsilon\\
(1,0)&\quad\text{with probability } \epsilon\\
(0,0)&\quad\text{with probability }\epsilon.\\
\endalign$$
This measure satisfies (1.3) if $4\epsilon\leq 1$, which we now assume.
Therefore $\mu_{\epsilon}S(t)$ is downward FKG for all $t\geq 0$ by assumption. 

We will use the shorthand $\nu(abc)$ to mean $\nu\{\eta:\eta(x)=a, \eta(y)=b,
\eta(z)=c\}$. The quantity
$$\mu_{\epsilon} S(t)(110)\mu_{\epsilon} S(t)(000)-\mu_{\epsilon} S(t)(100)\mu_{\epsilon}
S(t)(010)\tag 3.6$$ is zero at $t=0$, and is nonnegative for $t\geq 0$ since $\mu_{\epsilon}
S(t)$ is downward FKG. Therefore, its derivative is nonnegative at $t=0$. To write down this
dervative, let
$$\align\delta_{abc}(x)&=E_{\mu_{\epsilon}}[\delta(x,\cdot)\mid
\eta(x)=a,\eta(y)=b,\eta(z)=c],\\
\delta_{abc}(y)&=E_{\mu_{\epsilon}}[\delta(y,\cdot)\mid \eta(x)=a,\eta(y)=b,\eta(z)=c],\\
\delta_{abc}(z)&=E_{\mu_{\epsilon}}[\delta(z,\cdot)\mid \eta(x)=a,\eta(y)=b,\eta(z)=c],\\
\endalign$$
with $\beta_{abc}(x),\beta_{abc}(y)$ and $\beta_{abc}(z)$ defined similarly.
Then the derivative of (3.6) at $t=0$ is
$$\gather\mu(110)\bigg[\mu(100)\delta_{100}(x)+\mu(010)\delta_{010}(y)+\mu(001)\delta_{001}(z)-
\mu(000)[\beta_{000}(x)+\beta_{000}(y)+\beta_{000}(z)]\bigg]\\+
\mu(000)\bigg[\mu(010)\beta_{010}(x)+\mu(100)\beta_{100}(y)+\mu(111)\delta_{111}(z)-
\mu(110)[\delta_{110}(x)+\delta_{110}(y)+\beta_{110}(z)]\bigg]\\-
\mu(100)\bigg[\mu(110)\delta_{110}(x)+\mu(000)\beta_{000}(y)+\mu(011)\delta_{011}(z)-
\mu(010)[\beta_{010}(x)+\delta_{010}(y)+\beta_{010}(z)]\bigg]\\-
\mu(010)\bigg[\mu(000)\beta_{000}(x)+\mu(110)\delta_{110}(y)+\mu(101)\delta_{101}(z)-
\mu(100)[\delta_{100}(x)+\beta_{100}(y)+\beta_{100}(z)]\bigg],
\endgather$$
where we have omitted the subscript $\epsilon$.
Dividing this by $\epsilon$ and letting $\epsilon\downarrow 0$ yields
$$\delta_{111}(z)\geq\delta_{011}(z).$$
In other words,
$$\int_{\{\eta(x)=1\}}\delta(z,\eta)d\mu_0\geq\int_{\{\eta(x)=0\}}\delta(z,\eta)d\mu_0.$$
Since $\delta(z,\eta)$ is decreasing in $\eta$, it follows that
$\delta(z,\eta)=\delta(z,\eta_x)$ for all $\eta$ such that $\eta(x)=0,\eta(y)=1$.
Letting $x$ and $y$ vary, we see that $\delta(z,\eta)$ is constant on
$\{\eta:\eta\not\equiv 0\}$ as required.
\enddemo
\demo{Proof of Theorem 1.6(b)} Since (1.3) implies downward FKG, Proposition 3.5
can be applied to conclude that $\delta(z,\eta)$ is constant on
$\{\eta:\eta\not\equiv 0\}$. Interchanging the roles of 0's and 1's,
we see that $\beta(z,\eta)$ is constant on $\{\eta:\eta\not\equiv 1\}$.
This argument is correct, since the hypothes1s of Theorem 1.6(b) is symmetric
in 0's and 1's.

To complete the proof that $\delta(z,\eta)$ and $\beta(z,\eta)$ are
independent of $\eta$, we need to assume that $S$ has at least
four points. Let $\delta(z)$ be the
value of $\delta(z,\eta)$ on the set $\{\eta:\eta\not\equiv 0\}$, and 
let $\beta(z)$ be the value of $\beta(z,\eta)$ on the set $\{\eta:\eta\not\equiv 1\}$.
Take $\mu$ to be any product measure satisfying
$0<\mu\{\eta:\eta(z)=1\}<1$ for all $z\in S$.
Fix distinct sites $x,y$ and take nonempty sets $A,B$ so that
$A,B, \{x,y\}$ form a partition of $S$. Let $\eta$ and $\zeta$ be the
configurations
$$\eta(z)=\cases 0\quad&\text{ on }A\cup\{x\}\\1\quad&\text{ on }B\cup\{y\}\endcases$$
and
$$\zeta(z)=\cases 0\quad&\text{ on }B\cup\{x\}\\1\quad&\text{ on }A\cup\{y\}.\endcases$$
Then
$$\mu S(t)(\eta\vee\zeta)\mu S(t)(\eta\wedge\zeta)-\mu S(t)(\eta)\mu
S(t)(\zeta)\tag 3.7$$
is zero at $t=0$, and hence by assumption, its derivative must be nonnegative
at $t=0$. Since (1.3) holds with equality for any product measure,
and independent flip processes preserve the class of product measures,
the derivative of (3.7) would be zero if the spin system had
independent flips. This observation leads to a lot of
cancellation in the derivative of (3.7) in the present case,
since the process almost has independent flips. Letting 0 and 1 denote the
configurations that are identically 0 and 1 respectively, we find after
using this cancellation that
$$\gathered\frac d{dt}\bigg[\mu S(t)(\eta\vee\zeta)\mu S(t)(\eta\wedge\zeta)-\mu
S(t)(\eta)\mu S(t)(\zeta)\bigg]\bigg|_{t=0}\\=\mu(\eta\vee\zeta)\mu(\eta\wedge\zeta)
\bigg[\delta(y)-\delta(y,0)+\beta(x)-\beta(x,1)\bigg].
\endgathered\tag 3.8$$
Since the spin system is attractive (by Proposition 2.1), $\delta(y)\leq\delta(y,0)$ and
$\beta(x)\leq\beta(x,1)$, so it follows from the nonnegativity of (3.8) that 
$\delta(y)=\delta(y,0)$ and $\beta(x)=\beta(x,1)$ as required.
\enddemo
\heading 4. Preservation of the downward FKG property\endheading
This section is devoted to the proof of Theorem 1.8. We begin with a simple lemma.
The inequality in (4.2) below refers to stochastic monotonicity.

\proclaim {Lemma 4.1} Suppose $\mu$ is downward FKG and $A\subset B$. Then
$$\mu\{\cdot\mid\eta\equiv 0\text{ on }B\}\leq \mu\{\cdot\mid\eta\equiv 0\text{ on }A\},\tag
4.2$$
and any convex combination of these two conditional measures is associated.
\endproclaim

\demo{Proof}Let $f$
be an increasing function. By the downward FKG property, $f$ and $1_{\{\eta\equiv 0\text{ on
}B\}}$ are negatively correlated with respect to $\mu\{\cdot\mid\eta\equiv 0\text{ on }A\}$.
This gives
$$\int fd\mu\{\cdot\mid\eta\equiv 0\text{ on }B\}\leq \int fd\mu\{\cdot\mid\eta\equiv 0\text{
on }A\}$$
as required for (4.2). For the second statement, use Proposition 2.22 on page 83
of Liggett (1985).
\enddemo

\demo{Proof of Theorem 1.8(a)} By Proposition 1.14, it suffices to prove Theorem 1.8(a) for spin 
systems that have nonzero transition rates at only one site $u\in S$, and at that site,
either $\beta(u,\eta)\equiv 0$ and $\delta(u,\eta)\equiv 1$, or $\delta(u,\eta)\equiv 0$,
$\beta(x,\eta)\equiv 0$ for $x\neq u$, and $\beta(u,\eta)=1_{\{\eta\not\equiv 0\text
{ on }A\}}$ for a fixed $A\subset S$.

We begin then by considering the spin system with $\delta(u,\eta)\equiv 1$ for a fixed $u$, and
all other rates zero. Suppose that $\mu$ is downward FKG. If $u\notin A$, then the
evolution commutes with the operation of conditioning on $\{\eta\equiv 0\text{ on }A\}$, so that
$\mu S(t)\{
\cdot\mid\eta\equiv 0\text{ on }A\}$ is associated for all $t\geq 0$
by Theorem 1.4. So, we may assume that $u\in A$. If $\eta(u)=0$, then
$$\mu S(t)(\eta)=\mu(\eta)+\mu(\eta_u)(1-e^{-t}),$$
so that
$$\mu S(t)\{\cdot\mid\eta\equiv 0\text{ on }A\}=\lambda\mu\{\cdot\mid\eta\equiv 0\text
{ on }A\}+(1-\lambda)\mu\{\cdot\mid\eta\equiv 0\text{ on }A\backslash\{u\}\},$$
where
$$\lambda=\frac{e^{-t}\mu\{\eta\equiv 0\text{ on }A\}}
{e^{-t}\mu\{\eta\equiv 0\text{ on }A\}+(1-e^{-t})\mu\{\eta\equiv 0\text{ on }
A\backslash\{u\}\}}.$$
Therefore, $\mu S(t)\{\cdot\mid\eta\equiv 0\text{ on }A\}$ is associated
by Lemma 4.1.

Turning to the second case, assume now that all flip rates are zero, except that
$\beta(u,\eta)=1_{\{\eta\not\equiv 0\text{ on }A\}}$ for a particular $u\in S$ and
$A\subset S\backslash\{u\}$. Let $\mu$ be downward FKG. We need to check that the
measure $\mu S(t)\{\cdot\mid\eta\equiv 0\text{ on }B\}$ is associated for
every $B\subset S$. If $u\notin B$, this is a consequence of Theorem 1.4, since
conditioning on $\{\eta\equiv 0\text{ on }B\}$ commutes with $S(t)$. So, we
may assume that $u\in B$. In this case,
$$\int_{\{\eta\equiv 0\text{ on }B\}}fd\mu S(t)=\int_{\{\eta\equiv 0\text{ on }B\}}f(\eta)
e^{-t\beta(u,\eta)}d\mu $$
for every $f$. Writing 
$$e^{-t\beta(u,\eta)}=e^{-t}+(1-e^{-t})1_{\{\eta\equiv 0\text{ on }A\}},$$
it follows that
$$\mu S(t)\{\cdot\mid \eta\equiv 0\text{ on }B\}=\lambda
\mu\{\cdot\mid\eta\equiv 0\text{ on }B\}+(1-\lambda)
\mu\{\cdot\mid\eta\equiv 0\text{ on }A\cup B\},$$
where
$$\lambda=\frac{e^{-t}\mu\{\eta\equiv 0\text{ on }B\}}
{e^{-t}\mu\{\eta\equiv 0\text{ on }B\}+(1-e^{-t})\mu\{\eta\equiv 0\text{ on }A\cup B\}}.$$
Therefore $\mu S(t)\{\cdot\mid \eta\equiv 0\text{ on }B\}$ is associated
by Lemma 4.1.
\enddemo

\demo{Proof of Theorem 1.8(b)} Since (1.3) implies downward FKG, this
follows immediately from Proposition 3.5.
\enddemo

\demo{Proof of Theorem 1.8(c)} Since every product measure is downward FKG (by Theorem 1.2)
and every measure that is downward FKG is associated, the fact that
$\beta(u,\eta)$ is increasing in $\eta$ is a consequence of Proposition 2.1. So, it
remains to prove (1.10). It is sufficient to check it in the case that $\eta$ and $\zeta$
differ at only two sites; call them $x$ and $y$. Since $\beta(u,\gamma)$ does not depend
on $\gamma(u)$, we may assume that $x,y,u$ are distinct. 

Take $\mu$ to be a product measure. By assumption, $\mu S(t)$ is downward FKG for all 
$t\geq 0$. Therefore, using the definition of downward FKG with the conditioning on
$\{\eta(u)=0\}$,
$$\gather\mu S(t)\{\eta(x)=1,\eta(y)=1,\eta(u)=0\}\mu S(t)\{\eta(x)=0,\eta(y)=0,\eta(u)=0\}-
\\\mu S(t)\{\eta(x)=1,\eta(y)=0,\eta(u)=0\}\mu
S(t)\{\eta(x)=0,\eta(y)=1,\eta(u)=0\}\endgather$$
is nonnegative for all $t\geq 0$, and is zero at $t=0$. It follows that its derivative
is nonnegative at $t=0$. Writing this out, we see that
$$\gather E[\beta(u,\eta)\mid\eta(x)=1,\eta(y)=0]+E[\beta(u,\eta)\mid\eta(x)=0,\eta(y)=1]\geq\\
E[\beta(u,\eta)\mid\eta(x)=1,\eta(y)=1]+E[\beta(u,\eta)\mid\eta(x)=0,\eta(y)=0],\endgather$$
where the conditional expectations are with respect to $\mu$. Since $\mu$ is an
arbitrary product measure, we can conclude that
$$\beta(u,\gamma_x)+\beta(u,\gamma_y)\geq \beta(u,\gamma_{x,y})+\beta(u,\gamma)$$
for every $\gamma$ such that $\gamma(x)=\gamma(y)=0$. 
But this is exactly (1.10) with $\eta=\gamma_x$ and $\zeta=\gamma_y$.
\enddemo

\heading 5. Preservation of the DCA property\endheading
This section is devoted to the proof of Theorem 1.13. 
\demo{Proof of Theorem 1.13(a)} By Proposition 1.14, it suffices to consider a spin system
with nonzero flip rates only at one site $u$, and at that site, only $\delta(u,\eta)$
or $\beta(u,\eta)$ is not identically zero. We assume this, and also that the nonzero rates satisfy
(1.9) or (1.10) and attractiveness, in the two cases. 

The proof relies on three facts:

Fact I. If $f$ and $g$ are increasing and $h$ is positive, then
$$[S(t)h][S(t)(fgh)]\geq [S(t)(fh)][S(t)(gh)].\tag 5.1$$

Fact II. If $h$ is positive, decreasing and satisfies (1.11), then $S(t)h$ has the
same three properties.

Fact III. If $f$ is increasing, and $h$ is positive and satisfies (1.11), then
the function $f_t$ defined by
$$f_t(\eta)=\frac{S(t)(fh)(\eta)}{S(t)h(\eta)}$$
is increasing in $\eta$.
\medskip
\noindent{\bf Remark.} To see why the proof of Theorem 1.4 is easier than the proof of
Theorem 1.13(a), note that if $h\equiv 1$, then Facts II and III are immediate. (Fact I
is (2.20) on page 81 of Liggett (1985).)
\medskip
First we will deduce Theorem 1.13(a) from these three statements. Suppose $\mu$ is DCA,
$h$ is positive, decreasing and satisfies (1.11), and $f$ and $g$ are increasing.
Then $\mu_{S(t)h}$ is associated by Fact II. Applying this association to the increasing
functions $f_t$ and $g_t$, which are increasing by Fact III, we see that
$$\int S(t)hd\mu\int\frac{S(t)(fh)S(t)(gh)}{S(t)h}d\mu\geq\int S(t)(fh)d\mu\int S(t)(gh)d\mu.$$
Combining this with Fact I gives
$$\int S(t)hd\mu\int S(t)(fgh)d\mu\geq\int S(t)(fh)d\mu\int S(t)(gh)d\mu,$$
which can be rewritten as
$$\int hd[\mu S(t)]\int fghd[\mu S(t)]\geq\int fhd[\mu S(t)]\int ghd[\mu S(t)],$$
so that $\mu S(t)$ is DCA, as required.

We turn now to the proofs of Facts I and III. Since 

(a) monotonicity of $h$ is not assumed here,

(b) condition (1.11) is symmetric in 0's and 1's, 

\noindent and 

(c) the assumptions on the death
rates are more stringent than those on the birth rates, 

\noindent it suffices to prove this in
the case that the nonzero flip rates are $\beta(u,\eta)$. For a fixed $t$,
let $b(\eta)=e^{-t\beta(u,\eta)}.$ This is the probability that
$\eta_t(u)=0$ if the initial configuration satisfies $\eta(u)=0$. Then, for any function $f$,
$$S(t)f(\eta)=\cases f(\eta)&\text{ if }\eta(u)=1\\b(\eta)f(\eta)+(1-b(\eta))f(\eta_u)
&\text{ if }\eta(u)=0.\endcases$$
Therefore, the two sides of (5.1) are equal when the $\eta$ at which they are evaluated 
satisfies $\eta(u)=1$. If $\eta(u)=0$, then 
$$\align [S(t)h](\eta)&[S(t)(fgh)](\eta)- [S(t)(fh)](\eta)[S(t)(gh)](\eta)\\
&=\big[b(\eta)h(\eta)+(1-b(\eta))h(\eta_u)\big]\big[b(\eta)f(\eta)g(\eta)h(\eta)+(1-b(\eta))
f(\eta_u)g(\eta_u)h(\eta_u)\big]\\&\qquad-
\big[b(\eta)f(\eta)h(\eta)+(1-b(\eta))
f(\eta_u)h(\eta_u)\big]\big[b(\eta)g(\eta)h(\eta)+(1-b(\eta))
g(\eta_u)h(\eta_u)\big]\\&=b(\eta)(1-b(\eta))h(\eta)h(\eta_u)[f(\eta_u)-f(\eta)][g(\eta_u)
-g(\eta)]\geq 0.\endalign$$
This proves Fact I.

For Fact III, we must show that for any $v$, if $\eta(v)=0$, then
$$\frac{S(t)(fh)(\eta_v)}{S(t)h(\eta_v)}\geq \frac{S(t)(fh)(\eta)}{S(t)h(\eta)}.\tag 5.2$$
If $v=u$, the difference between the left and right sides of (5.2) is
$$f(\eta_u)-\frac{b(\eta)f(\eta)h(\eta)+(1-b(\eta))f(\eta_u)h(\eta_u)}
{b(\eta)h(\eta)+(1-b(\eta))h(\eta_u)}=\frac{b(\eta)h(\eta)[f(\eta_u)-f(\eta)]}
{b(\eta)h(\eta)+(1-b(\eta))h(\eta_u)}\geq 0.$$
If $v\neq u$ and $\eta(u)=1$, then the difference between the left and right sides of (5.2) is
$$f(\eta_v)-f(\eta)\geq 0.$$
If $v\neq u$ and $\eta(u)=0$, then the difference between the left and right sides of (5.2) is
$$\frac{b(\eta_v)f(\eta_v)h(\eta_v)+(1-b(\eta_v))f(\eta_{u,v})h(\eta_{u,v})}
{b(\eta_v)h(\eta_v)+(1-b(\eta_v))h(\eta_{u,v})}-\frac{b(\eta)f(\eta)h(\eta)+(1-b(\eta))f(\eta_u)h(\eta_u)}
{b(\eta)h(\eta)+(1-b(\eta))h(\eta_u)}.$$
Putting this over a common denominator, the resulting numerator is
$$\gathered
b(\eta)b(\eta_v)h(\eta)h(\eta_v)[f(\eta_v)-f(\eta)]+(1-b(\eta))(1-b(\eta_v))h(\eta_u)h(\eta_{u,v})[f(\eta_{u,v}
-f(\eta_u)]\\+b(\eta_v)(1-b(\eta))h(\eta_u)h(\eta_v)[f(\eta_v)-f(\eta_u)]+
b(\eta)(1-b(\eta_v))h(\eta)h(\eta_{u,v})[f(\eta_{u,v}-f(\eta)].
\endgathered\tag 5.3$$
Since $f$ is increasing, all terms but the third are nonnegative. To see that the nonnegative terms
compensate for the potentially negative one, we proceed as follows. The values of $f$ at the
four configurations that appear in (5.3) satisfy $f(\eta)\leq f(\eta_u), f(\eta_v)\leq f(\eta_{u,v}).$
Since they appear linearly in (5.3), it is enough to check the nonnegativity of (5.3) in case
$f(\eta)=0, f(\eta_{u,v})=1$, and $(f(\eta_u),f(\eta_v))$ takes one of the four values
$(0,0),(1,0),(0,1),(1,1)$. The summands in (5.3) are all nonnegative except in the case $f(\eta_u)=1,
f(\eta_v)=0.$ In this case, (5.3) becomes
$$b(\eta)(1-b(\eta_v))h(\eta)h(\eta_{u,v})-b(\eta_v)(1-b(\eta))h(\eta_u)h(\eta_v).$$
But, this is nonnegative since $h(\eta)h(\eta_{u,v})\geq h(\eta_u)h(\eta_v)$ (by (1.11))
and $b(\eta)\geq b(\eta_v)$ (since the birth rates $\beta(u,\eta)$ are increasing).

It remains to prove Fact II. We will again carry out the proof in the case that the nonzero
flip rates are the birth rates at site $u$. To check the result for the case of nonzero
death rates, simply note that the monotonicity of $h$ is not needed in the proof that
$S(t)h$ satisfies (1.11) if $\beta(u,\eta)$ is constant, so that the interchange of
roles of zeros and ones can be used again. The fact that $S(t)h$ is positive is
clear, and the fact that it is decreasing follows from the monotonicity of $\beta(u,\eta)$.

To verify that $S(t)h$ satisfies (1.11), it is enough to check that if $\eta$ satisfies
$\eta(v)=\eta(w)=0$ for two distinct sites $v$ and $w$, then
$$S(t)h(\eta_{v,w})S(t)h(\eta)\geq S(t)h(\eta_v)S(t)h(\eta_w).\tag 5.4$$
If $w=u$, then the difference between the left and right sides of (5.4) is
$$b(\eta_v)[h(\eta)h(\eta_{u,v})-h(\eta_u)h(\eta_v)]+[b(\eta)-b(\eta_v)]h(\eta_{u,v})[h(\eta)-h(\eta_u)],$$
which is nonnegative if $b$ is decreasing (i.e., $\beta(u,\eta)$ is increasing), $h$ is
decreasing and satisfies (1.11). Note that the monotonicity of $h$ is not needed if $b$
is constant. Now suppose that $u,v,w$ are all distinct. Again, we may assume that $\eta(u)=0$,
since otherwise (5.4) is automatic. Then the difference between the left and right sides
of (5.4) is
$$\gathered\big[b(\eta_{v,w})h(\eta_{v,w})+(1-b(\eta_{v,w}))h(\eta_{u,v,w})\big]
\big[b(\eta)h(\eta)+(1-b(\eta))h(\eta_{u})\big]\\
-\big[b(\eta_{v})h(\eta_{v})+(1-b(\eta_{v}))h(\eta_{u,v})\big]
\big[b(\eta_w)h(\eta_w)+(1-b(\eta_w))h(\eta_{u,w})\big].\endgathered\tag 5.5$$
We need to check that this is nonnegative whenever
$$0<b(\eta_{v,w})\leq b(\eta_v),b(\eta_w)\leq b(\eta)\leq 1, \text{ and }
b(\eta_v)b(\eta_w)\leq b(\eta)b(\eta_{v,w}).\tag 5.6$$

The easiest case to consider is (i), in which  $\{b(\eta_v),b(\eta_w)\}= \{b(\eta),b(\eta_{v,w})\}$.
In this case, we can define a probability measure $\sigma$ on $\{0,1\}^3$ by
$$\gather
\sigma(111)=(1-b(\eta_{v,w}))h(\eta_{u,v,w}),\quad \sigma(110)=b(\eta_{v,w})h(\eta_{v,w})\\
\sigma(101)=(1-b(\eta_{v}))h(\eta_{u,v}),\quad\sigma(100)=b(\eta_{v})h(\eta_{v})\\
\sigma(011)=(1-b(\eta_{w}))h(\eta_{u,w}),\quad\sigma(010)=b(\eta_{w})h(\eta_{w})\\
\sigma(001)=(1-b(\eta))h(\eta_{u}),\quad\sigma(000)=b(\eta)h(\eta),
\endgather$$
and then normalizing it to sum to 1. This measure satisfies the FKG lattice condition, and
hence by Theorem 1.2, is associated. But in this case, (5.5) is just a constant multiple
of the covariance of the first two coordinates relative to $\sigma$, so it is nonnegative.
Note that in this case, we have not used the monotonicity of $h$.

The next case is (ii), in which $b(\eta_v)=b(\eta_w)=b(\eta_{v,w})$ -- call this
common value $a$. 
Since (5.5) is linear in $b(\eta)$, it suffices to check its nonnegativity in
the extreme cases $b(\eta)=1$, and $b(\eta)=a$.
The latter case is a special case of case (i). So, we may assume that $b(\eta)=1$.
In this case, (5.5) is a quadratic polynomial in $a$, and the coefficient of
$a^2$ is
$$-[h(\eta_v)-h(\eta_{u,v})][h(\eta_w)-h(\eta_{u,w})].$$
This is nonpositive by the monotonicity of $h$. (Note that $h$ monotone in either
direction would be enough here.) So, it is enough to check the nonnegativity at the
two extreme cases, $a=0$ and $a=1$. The case $a=1$ is again a special case of case (i).
If $a=0$, then (5.5) becomes
$$h(\eta_{u,v,w})h(\eta)-h(\eta_{u,v})h(\eta_{u,w}),$$
which is nonnegative by (1.11) and the fact that $h$ is {\it decreasing}. This completes
the consideration of case (ii).

Now think of (5.5) as a function of the variables $x=b(\eta_v)$ and $y=b(\eta_w)$ for fixed values of
$0<b(\eta_{v,w})<b(\eta)\leq 1$. Since (5.5) is bilinear in these two variables, it
suffices to check its nonnegativity on the boundary of the region defined by (5.6).
We have already checked it at three points on the boundary, in cases (i) and (ii) above. 
The boundary consists of two line segments and one curve, so again by bilinearity, 
the nonnegativity follows on the line segments from its nonnegativity at the endpoints
of those line segments. It suffices then to consider the case in which $xy=b(\eta)b(\eta_{v,w})=A$.
Replacing $y$ by $A/x$ in (5.5) and expanding, we see that the only dependence on $x$ is
in the terms
$$-xh(\eta_{u,w})[h(\eta_v)-h(\eta_{u,v})]-\frac Axh(\eta_{u,v})[h(\eta_w)-h(\eta_{u,v})].$$
Since $h$ is decreasing, this is a concave function of $x$, so (5.5) will be proved to
be nonnegative once it is nonnegative at the endpoints of the interval of $x$'s that
are relevant. But this again corresponds to case (i), so the proof is complete.
\enddemo

\demo{Proof of Theorem 1.13(b)} This is again a consequence of Proposition 3.5, 
since
$$\mu \text{ satisfies (1.3)}\Rightarrow \mu\text{ DCA }\Rightarrow\mu S(t) \text{ DCA }
\Rightarrow\mu S(t)\text{ downward FKG}.$$
\enddemo
The proof of Theorem 1.13(c) is the same as that of Theorem 1.8(c).

\heading 6. The case of three sites\endheading
In this section, we take $S=\{1,2,3\}$, and find necessary and sufficient conditions  
for $\mu$ to satisfy each of the four properties appearing in (1.12). It will follow
that for three sites, the first and third implications in (1.12) are strict, while the
second is an equivalence. 

We will use the following notation in this section:
$$\gather a=\mu(111),\\b_1=\mu(011),\quad b_2=\mu(101),\quad b_3=\mu(110),\\
c_1=\mu(100),\quad c_2=\mu(010),\quad c_3=\mu(001),\\d=\mu(000).
\endgather$$
To state the necessary and sufficient conditions for the four properties of interest,
consider the following sets of inequalities:
$$\gathered a(c_2+c_3+d)\geq b_1(b_2+b_3+c_1),\\a(c_1+c_3+d)\geq b_2(b_1+b_3+c_2),\\
a(c_1+c_2+d)\geq b_3(b_1+b_2+c_3),\endgathered\tag A$$
$$\gathered d(b_2+b_3+a)\geq c_1(c_2+c_3+b_1),\\d(b_1+b_3+a)\geq c_2(c_1+c_3+b_2),\\
d(b_1+b_2+a)\geq c_3(c_1+c_2+b_3),\endgathered\tag B$$
$$\gathered(b_1+a)(c_1+d)\geq(c_3+b_2)(b_3+c_2),\\(b_2+a)(c_2+d)\geq(c_1+b_3)(b_1+c_3),
\\(b_3+a)(c_3+d)\geq(c_2+b_1)(b_2+c_1),\endgathered\tag C$$
$$b_1 d\geq c_2 c_3,\qquad  b_2 d\geq c_1 c_3,\qquad b_3 d\geq c_1 c_2,\tag D$$
and
$$c_1 a\geq b_2 b_3,\qquad  c_2 a\geq b_1 b_3,\qquad c_3 a\geq b_1 b_2.\tag E$$

Here are the necessary and sufficient conditions:

\proclaim{Proposition 6.1} (a) $\mu$ satisfies the FKG lattice condition if and only if (D)
and (E) hold.

(b) $\mu$ satisfies the DCA property if and only if (A), (C) and (D) hold.

(c) $\mu$ satisfies the downward FKG property if and only if (A), (C) and (D) hold.

(d) $\mu$ is associated if and only if (A), (B) and (C) hold.
\endproclaim

\noindent {\bf Remark.} For an example to show that the first implication in (1.12) is strict, take
$\epsilon$ small and
$$a=b_1=b_2=b_3=1/6,\quad c_1=c_2=c_3=\epsilon,\quad d=1/3.$$
To show that the third implication in (1.12) is strict, take
$$a=1/3,\quad b_1=b_2=b_3=\epsilon,\quad c_1=c_2=c_3=d=1/6.$$

The proof of part (a) of the theorem is immediate. We turn now to the other parts.

\demo{Proof of Proposition 6.1 (d)} 
To check that association implies (A), (B) and (C), let
$f_i(\eta)=\eta(i), g_i(\eta)=\prod_{j\neq i}\eta(j),$ and $h_i(\eta)=$  the
indicator of the event $\{\sum_{j\neq i}\eta(j)\geq 1\}$. Using $cov$ to denote the covariance
with respect to $\mu$, we have the following:
$$\gathered cov(f_1,g_1)=a(c_2+c_3+d)-b_1(b_2+b_3+c_1),\\cov(f_1,h_1)=d(b_2+b_3+a)-c_1(c_2+c_3+b_1),\\
cov(f_2,f_3)=(b_1+a)(c_1+d)-(c_3+b_2)(b_3+c_2).
\endgathered\tag 6.2$$
If $\mu$ is associated, each of these covariances is nonnegative.
This gives the first inequality in each of (A), (B) and (C). The others are obtained by
permuting the coordinates.

The proof of the converse is longer. Assume that (A), (B) and (C) hold. 
By (6.2) and the corresponding inequalities obtained by permuting the
coordinates, we may assume that
$$cov(f_i,g_i)\geq 0,\quad cov(f_i,h_i)\geq 0,\quad\text{and}\quad cov(f_i,f_j)\geq 0\tag 6.3$$
for all choices of $i$ and $j$. We need
to check that $cov(f,g)\geq 0$ for all increasing functions $f$ and $g$. It is
sufficient to check this when $f$ and $g$ are both increasing indicator
functions, since any increasing function can be written as a positive
linear combination of increasing indicator functions. Writing 
$$cov(f,g)=\frac 12\sum_{\eta,\zeta}[f(\eta)-f(\zeta)][g(\eta)-g(\zeta)]\mu(\eta)\mu(\zeta),$$
we see that $cov(f,g)\geq 0$ if $[f(\eta)-f(\zeta)][g(\eta)-g(\zeta)]\geq 0$
for all $\eta,\zeta$. Since we are assuming that $f$ and $g$ are indicators, 
these products can only take the values $-1,0,+1$. Therefore, we may assume that
at least one of these products is $-1$. The product is nonnegative
whenever $\eta$ and $\zeta$ are comparable. Therefore, by permuting coordinates and/or interchanging
the roles of 0's and 1's, we see that there are two cases to consider:
\medskip
\noindent {\it Case 1}. $[f(101)-f(011)][g(101)-g(011)]=-1$, in which case we may
assume 
$$f(101)=1,\quad f(011)=0,\quad g(101)=0,\quad g(011)=1.$$
The monotonicity of $f$ and $g$ forces $f$ and $g$ to also take the
following values:
$$f(111)=g(111)=1,\quad\text{and}\quad f(010)=f(001)=f(000)=g(100)=g(001)=g(000)=0.$$
 So, there are only three possibilities for $f$ and $g$:
$$f=g_2, f=f_1h_1,\text{ or }f=f_1\text{ and }g=g_1,g=f_2 h_2,\text{ or }g_2=f_2.$$
The corresponding nine covariances are nonnegative by (6.3), since
$$\gathered cov(g_i,g_j)\geq cov(g_i,f_i),\quad cov(g_i,f_ih_i)\geq cov(g_i,f_i)\\
cov(f_ih_i,f_jh_j)\geq cov (f_i,f_j),\quad cov (f_ih_i,f_j)\geq cov(f_i,f_j)\endgathered\tag 6.4$$
for $i\neq j$. These all follow from the easy fact that if $G_i$ and $H_i$ are
sets satisfying $G_i\subset H_i$ for $i=1,2$ and $G_1\cap G_2=H_1\cap H_2$, then
$cov(1_{G_1},1_{G_2})\geq cov(1_{H_1},1_{H_2}).$
\medskip
\noindent {\it Case 2}. $[f(011)-f(100)][g(011)-g(100)]=-1$, in which case we may
assume 
$$f(011)=1,\quad f(100)=0,\quad g(011)=0,\quad g(100)=1.$$
Using the monotonicity of $f$ and $g$ as before, it follows that $g=f_1$ and $f$ is one of the
following:
$$g_1,\quad f_3 h_3,\quad f_3,\quad f_2 h_2,\quad h_1h_2h_3,\quad h_1h_2,\quad f_2,\quad h_1 h_3,\quad h_1.$$
The nonnegativity of the covariance of each of these with $f_1$ follows from (6.3), (6.4),
$$cov(h_1h_2h_3,f_i)\geq cov (h_i,f_i)\quad\text{and}\quad cov(h_i h_j,f_i)\geq cov(h_i,f_i)$$
for $i\neq j$.
\enddemo

\demo{Proof of Proposition 6.1 (c)} The statement that downward FKG implies (A), (C), and (D)
is an immediate consequence of part (d). So, assume now that (A), (C) and (D) hold. Since (D) holds,
it will be sufficient to show that $\mu$ is associated. And to do this, it suffices by part (d)
to show that (B) holds. To do so, multiply the first inequality in (A) by $d$ and then use
the second and third inequalities in (D) to conclude that
$$ad(c_2+c_3+d)\geq db_1(b_2+b_3+c_1)\geq b_1 c_1(c_2+c_3+d),$$
and therefore that $ad\geq b_1c_1$. For future reference, we record the fact that
$$\text{(A) and (D) imply } ad\geq b_ic_i\text{ for }i=1,2,3.\tag 6.5$$
Now the first inequality in (B), for example, follows from (6.5) with $i=1$, together with
the last two inequalities in (D). 
\enddemo
\demo{Proof of Proposition 6.1 (b)} Now let $a,b_i,c_i,d$ be the probabilities of the various
configurations for the measure $\mu$, and $a^*,b_i^*,c_i^*, d^*$ be the corresponding probabilities
for the measure $\mu_h$, where $h$ is positive, decreasing, and satisfies (1.11). We need to
show that if the unstarred quantities satisfy (A), (C) and (D), then the starred quantities
satisfy (A), (B) and (C). By (1.11), the starred quantities satisfy (D), since the unstarred
quantities do. Therefore, once we have shown that the starred quantities satisfy (A) and (C)
they will automatically satisfy (B) by part (c) of the proposition.

We will now check that the starred quantities satisfy the first inequality in (A). To do so, 
define
$$\gather x_1=\frac{h(001)}{h(011)},\quad x_2=\frac{h(010)}{h(011)},\quad x_3=\frac{h(000)}{h(011)},\\
y_1=\frac{h(101)}{h(111)},\quad y_2=\frac{h(110)}{h(111)},\quad y_3=\frac{h(100)}{h(111)}.\endgather$$
The starred version of the first inequality in (A) then becomes
$$a(x_2 c_2+x_1 c_3+x_3 d)\geq b_1(y_1 b_2+y_2 b_3+y_3 c_1).\tag 6.6$$
Since $h$ is decreasing, we have
$$x_1,x_2,x_3,y_1,y_2,y_3\geq 1.\tag 6.7$$
Since $h$ satisfies (1.11), these quantities satisfy the following inequalities:
$$x_1\geq y_1,x_2\geq y_2,x_3\geq y_3,\quad y_3\geq y_1 y_2,x_3\geq x_1x_2,\quad y_1x_3\geq x_1y_3,y_2x_3\geq
x_2y_3.\tag 6.8$$
By the first three inequalities in (6.8), (6.6) will be a consequence of
$$a(y_2 c_2+y_1 c_3+y_3 d)\geq b_1(y_1 b_2+y_2 b_3+y_3 c_1).\tag 6.9$$
By (6.5), $ad\geq b_1c_1$. Therefore, by the fourth inequality in (6.8), the worst case
of (6.9) is that in which $y_3=y_1y_2.$ Consider this case and divide (6.9) by $y_1y_2$.
By (6.7), we see that we need to show that
$$a(z_1c_2+z_2c_3+d)\geq b_1(z_2b_2+z_1b_3+c_1)\tag 6.10$$
for all $0\leq z_1,z_2\leq 1$. Since the expressions on both sides of (6.10) are linear
in $z_1,z_2$, it suffices to check that (6.10) holds at the four corners of this square. If $z_1=z_2=0$,
(6.10) is true by (6.5). If $z_1=z_2=1$, (6.10) is just the first inequality in (A). 
The other two cases are $a(c_2+d)\geq b_1(b_3+c_1)$ and $a(c_3+d)\geq b_1(b_2+c_1)$. One is
obtained from the other by permutation of coordinates, so we will prove the first of these.

To do so, note that
$$\aligned [(b_3+c_1)(c_2+c_3+d)-&(b_2+b_3+c_1)(c_2+d)]\\&+[b_1(b_3+c_1)-(c_1+b_3)(b_1+c_3)+b_2(c_2+d)]=0.
\endaligned\tag 6.11$$
Therefore one of the expressions in brackets is nonpositive. If it is the first one, then
$$(b_3+c_1)(c_2+c_3+d)\leq (b_2+b_3+c_1)(c_2+d),$$
which when combined with the first inequality in (A) implies $a(c_2+d)\geq b_1(b_3+c_1)$.
If it is the second expression in brackets in (6.11) that is nonpositive, then
$$b_1(b_3+c_1)\leq(c_1+b_3)(b_1+c_3)-b_2(c_2+d),$$
which when combined with the second inequality in (C) implies $a(c_2+d)\geq b_1(b_3+c_1)$.
Therefore, this last inequality is true in either case.

The starred version of the third inequality in (C) (we choose the third one to check instead
of the first so that we can use the same $x_i$'s and $y_i$'s as in the previous
argument) is
$$(b_3 y_2+a)(c_3 x_1+dx_3)\geq (c_2 x_2+b_1)(b_2y_1+c_1y_3).\tag 6.12$$
To check this inequality, start by writing
$$\align (b_3 y_2+a)(c_3 x_1+dx_3)-& (c_2 x_2+b_1)(b_2y_1+c_1y_3)\\
=&(b_3 y_2+a)(c_3 y_1+dy_3)- (c_2 y_2+b_1)(b_2y_1+c_1y_3)\\
&+c_3(b_3y_2+a)(x_1-y_1)+c_2(b_2y_1+c_1y_3)\bigg[\frac{y_2 x_3-x_2y_3}{y_3}\bigg]\\
&+\bigg[(db_3-c_1c_2)y_2+da-c_2b_2\frac{y_1y_2}{y_3}\bigg](x_3-y_3).
\endalign$$
Therefore, by (D), (6.5) and (6.8), it suffices to prove
$$(b_3 y_2+a)(c_3 y_1+dy_3)\geq (c_2 y_2+b_1)(b_2y_1+c_1y_3).\tag 6.13$$
By (D) and (6.5), the coefficient of $y_3$ on the left side of (6.13)
is at least as large as the coefficient of $y_3$ on the right side of (6.13).
Therefore, by (6.8), it suffices to consider the case $y_3=y_1y_2$.
Cancelling a common factor of $y_1$, we see by (6.7) that it suffices to check
$$(b_3 y_2+a)(c_3 +dy_2)- (c_2 y_2+b_1)(b_2+c_1y_2)\geq 0\tag 6.14$$
for $y_2\geq 1$. Write this polynomial as $p(y_2)$. Then $p(1)\geq 0$ by (C).
The coefficient of $y^2$ in $p(y)$ is $b_3d-c_2c_1$, which is nonnegative by (D).
So, it suffices to check that $p'(1)\geq 0$. But this follows from (C) and (D), since
$$\align d(c_3+d)p'(1)=&d^2[(b_3+a)(c_3+d)-(c_2+b_1)(b_2+c_1)]\\
&+(c_3+d)^2(b_3d-c_1c_2)+(b_1 d-c_2c_3)(b_2d-c_1c_3).\endalign$$
\enddemo
\bigskip
\centerline{\bf References}
\medskip

\ref \by V. Belitsky, P. Ferrari, N. Konno and T. M. Liggett
\paper A strong correlation inequality for contact processes and oriented
percolation\jour Stoch. Proc. Appl.\vol 67\yr 1997\pages 213--225\endref

\ref\by J. van den Berg, O. H\"aggstr\"om and J. Kahn\paper Some conditional
correlation inequalities for percolation and related processes
\yr 2005a\jour Rand. Structures Algorithms\endref

\ref\by J. van den Berg, O. H\"aggstr\"om and J. Kahn\paper Proof of a conjecture
of N. Konno for the 1D contact process
\yr 2005b\endref

\ref\by S. N. Ethier and T. G. Kurtz\book Markov Processes: Characterization
and Convergence\yr 1986\publ Wiley\endref

\ref\by P. C. Fishburn, P. G. Doyle and L. A. Shepp\paper The match set of
a random permutation has the FKG property\jour Ann. Probab.\yr 1988\vol 16
\pages 1194--1214\endref

\ref\by T. M. Liggett\book Interacting Particle Systems\publ Springer
\yr 1985\finalinfo reprinted in the series ``Classics
in Mathematics" in 2005\endref

\ref \by T. M. Liggett \paper Survival and coexistence in interacting particle systems
\inbook Probability and Phase Transition \publ Kluwer \yr 1994\pages 209--226\endref

\ref\by T. M. Liggett and J. E. Steif\paper Stochastic domination:
The contact process, Ising models and FKG measures\jour Ann. Inst. H. Poincar\'e Probab. Statist.\yr 2005\endref

\bigskip

\noindent Thomas M. Liggett

\noindent Department of Mathematics

\noindent University of California

\noindent Los Angeles CA 90095
\bigskip
\noindent Email: tml\@math.ucla.edu

\end